\documentclass[12pt]{article}

\usepackage{amsfonts, amssymb, amsmath, amsthm, epsf, cite}

\usepackage{fullpage}

\theoremstyle{plain}
\newtheorem{theorem}{Theorem}[section]
\newtheorem{lemma}{Lemma}[section]

\theoremstyle{definition}
\newtheorem{definition}{Definition}[section]

\newtheorem{remark}{Remark}[section]

\numberwithin{equation}{section} \numberwithin{figure}{section}

\renewcommand{\phi}{{\varphi}}

\newcommand{\cA}{{\mathcal A}}
\newcommand{\cV}{{\mathcal V}}

\newcommand{\cP}{{\mathcal P}}

\newcommand{\cR}{{\mathcal R}}

\newcommand{\cS}{{\mathcal S}}

\newcommand{\cW}{{\mathcal W}}

\newcommand{\cU}{{\mathcal U}}

\newcommand{\oQ}{{\overline Q}}
\newcommand{\oa}{{\overline a}}

\newcommand{\ou}{{\overline u}}
\newcommand{\ox}{{\overline x}}
\newcommand{\ow}{{\overline w}}
\newcommand{\ov}{{\overline v}}

\newcommand{\ux}{{\underline  x}}

\newcommand{\bP}{{\mathbf P}}

\newcommand{\bS}{{\mathbf S}}

\newcommand{\bu}{{\mathbf u}}

\newcommand{\bbR}{{\mathbb R}}
\newcommand{\bbN}{{\mathbb N}}

\let\phi=\varphi

\newcommand {\Var}{{\mathop{\rm Var}}}

\title{Well-posedness of parabolic equations containing hysteresis
with diffusive thresholds}

\author{Pavel Gurevich\footnote{Free University Berlin, Germany; Peoples' Friendship University of Russia, Russia; email: gurevichp@gmail.com},
Dmitrii Rachinskii\footnote{Department of Mathematical Sciences,
University of Texas at Dallas, USA \& Department of Applied
Mathematics, University College Cork, Ireland; email:
Dmitry.Rachinskiy@utdallas.edu}}

\begin{document}

\maketitle


\begin{abstract}
We study complex systems arising, in particular, in
population dynamics, developmental biology, and bacterial
metabolic processes, in which each individual element obeys a
relatively simple hysteresis law (a non-ideal relay). Assuming
that hysteresis thresholds fluctuate, we consider the arising
reaction-diffusion system. In this case, the spatial variable
corresponds to the hysteresis threshold. We describe the
collective behavior of such a system in terms of the Preisach
operator with time-dependent measure which is a part of the
solution for the whole system. We prove the well-posedness of the
system and discuss the long-term behavior of solutions.
\end{abstract}

\section{Introduction}

\begin{figure}[ht]
{\ \hfill\epsfxsize60mm\epsfbox{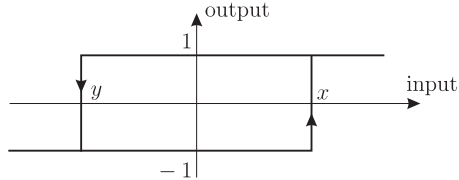}\hfill\ }
\caption{Non-ideal relay.}\label{figure1}
\end{figure}

In the paper, we develop a theory of reaction-diffusion equations
containing discontinuous hysteresis operator --- the so-called
non-ideal relay.  The non-ideal relay (or a bi-stable relay, or
lazy switch) is the most basic, yet non-trivial {\em hysteresis
operator}.   The state (output) of the non-ideal relay switches
from $-1$ to $1$ when the input   exceeds a threshold value
$x\in\mathbb R$ and switches back to state $-1$ when the input
drops below a smaller threshold value $y\in\mathbb R$. Hence,
variation of the input causes switching between two available
states in such a way that the current state depends on the input
history as long as the current value of the input falls within the
input's bi-stability range $(y,x)$. This behavior is illustrated
in Fig.~\ref{figure1}.

In particular, reaction-diffusion equations with non-ideal relay
arise in population dynamics, developmental biology (cell
differentiation processes), bacterial metabolic processes, etc.
The general feature of all those models is a hysteretic
interaction between several ``substances''. Depending on the
physical background, the substances are the amount of individuals,
concentration of proteins, density of nutrients, etc.

In  most existing models, the diffusion term in the
reaction-diffusion equation  appears because of a spatial motion
of the substances. In this paper, we suggest a principally new
mechanism of diffusion based on sporadical changes of the
thresholds of hysteresis operators. As we will explain, this
mechanism is not only natural, but also leads to a new dynamical
phenomena in the systems under consideration, e.g., appearance of
sign changing patterns for the states of non-ideal relays.

Let us  illustrate the above mentioned mechanism with the
following prototype example. Suppose we have a population of
bacteria with the environment being a mixture of two types of
nutrients --- lactose and glucose. Each bacterium, at a given
moment, can consume either lactose or glucose. In the first case,
we say that it is in the state $1$ and in the second case in the
state $-1$. For each bacterium at a given moment, there are two
(symmetric for simplicity) thresholds: $-x$ and~$x$, where $x>0$.
If the deviation of the relative concentration of lactose from the
value 1/2 in the mixture of the two nutrients is greater than $x$,
then the bacterium consumes lactose. If the deviation is less than
$-x$, it consumes glucose. If the deviation is between $-x$ and
$x$, then the state of the bacterium is defined according to the
discontinuous hysteresis law (see Fig.~\ref{figure1}, where $y=-x$
and the above deviation plays a role of input, while the state of
the bacterium plays a role of output).

Now the key (and quite natural) assumption is that each bacterium
can {\em sporadically change} its switching threshold or produce
offsprings with different thresholds. Assume that these changes
occur according to the Gaussian distribution centered at~$x$. Then
the density $u(x,t)$ of the biomass of bacteria with given
switching thresholds $\pm x$ at a moment $t$ will satisfy the
diffusion equation. We stress that the diffusion term appears not
due to a spatial motion of bacterium, but rather due to the {\em
diffusion of switching thresholds.}

Note that the collective impact of the whole population of
bacteria upon the environment (concentration of lactose and
glucose in our example) is mathematically described as an integral
of all discontinuous hysteresis  operators weighted with the
density $u(x,t)$. This integral can be interpreted as the Preisach
operator \cite{alexei} with a time dependent density $u(x,t)$. In
the classical Preisach operator, the density is usually time
independent and is assumed to be given. The typical difficulty in
this case is how to identify it in particular applications. Our
approach allows one to overcome this difficulty (at least in some
cases) because we treat the density as a component of the solution
itself. This gives rise to the question of the limiting behavior
of the density $u(x,t)$. For example, if one can describe a global
attractor, there is no need to (precisely) identify  the initial
density $u(x,0)$, since, after some transition period of time, it
will be in a vicinity of the known global attractor.

In the  present paper, we introduce, seemingly for the first time,
a notion of hysteresis (discontinuous non-ideal relay) with
diffusive thresholds. This leads us to an initial boundary-value
problem for a reaction-diffusion system including, as reaction
terms, discontinuous hysteresis relay operators and the integral
of those. We prove the well-posedness of this problem and
partially address an issue of the large time behavior of its
solutions.

The paper is organized as follows. In Sec.~\ref{secMotivation}, we
discuss various natural mechanisms leading to multi-stability in
biological systems, which can be mathematically described in terms
of hysteresis operators\footnote{We include this section for the
convenience of a reader interested in biological background of
problems with hysteresis. All the rest sections in the paper can
be read independently of Sec.~\ref{secMotivation}.}.

In Sec.~\ref{secModelDescription}, based on the above prototype
example of two-phenotype bacteria, we give a  rigorous model
formulation --- the reaction-diffusion system with discontinuous
hysteresis operators. The unknown functions are the density
$u(x,t)$ of the biomass of bacteria with given switching
thresholds $\pm x$ at a moment $t$, and the amounts $f_1(t)$ and
$f_{-1}(t)$ of the two types of nutrients. In the same section, we
rewrite the system in terms of the new unknown functions: the
density $u$, the total amount $v=f_1+f_{-1}$ of nutrients, and the
deviation $w=f_1/(f_1+f_{-1})-1/2$ of the relative concentration
of one type of nutrients from the value~1/2 in the mixture of the
two nutrients. Though our model is based on a particular
biological system, it accounts for a number of quite general
phenomena and can be adapted for other systems (cf.
Sec.~\ref{secMotivation}).

Sections~\ref{secStatesPreisach}--\ref{secLargeTimeBehavior} are
devoted to the analysis of the model. In
Sec.~\ref{secStatesPreisach}, we prove continuity properties of
hysteresis operators provided that the input function is
continuous in time and has bounded variation. Although each
non-ideal relay (which we call $\cR^x(\cdot)$ for fixed
thresholds~$\pm x$) is a discontinuous operator, it turns out that
an infinite collection of the relays with different
thresholds~$\pm x$ (which we call $\cR(\cdot)$) can be viewed as a
continuous operator with values in suitable $L_q$ (with respect to
$x$) spaces. An important issue, however, is that such an operator
is Lipschitz continuous for $q=1$ only, and it is H\"older
continuous for $q>1$. Similarly, the Preisach operator (which we
call $\cP(\cdot)$), i.e., the integral of the discontinuous relay
operators (weighted with the time dependent density function) over
all available thresholds $\pm x$ turns out to be a continuous
operator.

In Sec.~\ref{secWellPosedness}, we prove the well-posedness of the
model. The existence of solutions is proved via the Schauder
fixed-point theorem, where the continuity of the ``collective''
hysteresis $\cR$ in $L_q$ spaces  with $q>1$ is exploited. As we
said, the collective hysteresis $\cR$ is not Lipschitz continuous
for $q>1$. Therefore, the contraction mapping principle does not
apply, and the uniqueness should be proved separately. This is
done via the semigroup approach and additional estimates in the
$L_1$ space (where $\cR$ is Lipschitz continuous).

In Sec.~\ref{secLargeTimeBehavior}, we study the behavior of the
density $u(x,t)$ and the total amount of nutrients $v(t)$ as
$t\to\infty$. In particular, we prove that $v(t)$ monotonically
decreases and tends to $0$, while $u(x,t)$ converges (uniformly in
$x$) to a spatially homogeneous equilibrium. The large time
behavior of $w(t)$ as well as the limiting distribution of two
phenotypes $\lim\limits_{t\to\infty}{\mathcal R}^x(w)(t)$ remains
an open question.

The paper ends with Sec.~\ref{secDiscussion} containing some
discussion of the results and open questions as well as possible
modifications and extensions of the model.

\section{Motivation and biological
background}\label{secMotivation} The idea that epigenetic
differences such as those arising in the process of cell
differentiation can be attributed to multi-stability or
multi-stationarity of living forms seems to have been first
articulated by Max Delbr\"uk \cite{delb}. A classical example of
multi-stability in biology is multi-stable behavior of {\em
lac-operon} in {\em E-Coli}. {\em Lac-operon} is a collection of
genes associated with transport and metabolism of lactose in the
bacterium. Expression of these genes can be turned on by certain
small molecules that have been called inducers. Novick and Weiner
\cite{novi} as well as Cohn and Horibata \cite{cohn,cohn1},
relying on prior work of others \cite{mono,benz,spie}, effectively
demonstrated that two phenotypes each associated with ``on'' and
``off'' state of {\em lac-operon} expression can be obtained from
the same culture depending on the history of exposure to the
inducer. Moreover, both phenotypes remain stable through multiple
generations of the bacterial culture after the extracellular
concentration of the inducer is reduced to lower levels, but not
removed completely. This behavior resembles the definition of the
{\em non-ideal relay}
illustrated in Fig.~\ref{figure1}.

The above mentioned early findings on the hysteresis of the {\em
lac-operon} enzymes were consistent with earlier findings on
regulation of other enzymes in yeast \cite{wing}.  Recent
experiments using molecular biology methods (such as those
incorporating green fluorescent protein expression under the {\em
lac-operon} promoter) permitted to confirm and further study the
region of bi-stability of the {\em lac-operon} when multiple input
variables (TMG that acts as the inducer and glucose) are used to
switch the {\em lac-operon} genes on and off. Multi-stable gene
expression and hysteresis has been well-documented in a number of
natural as well as artificially constructed systems
\cite{3,5,10,15,17}.

Several important issues require further discussion when it comes to hysteresis and multi-stability in biological systems. First is the fact that reproductive rates in different phenotype states are frequently different in a given environment.  For example, the growth rate of one phenotype is high in a lactose rich environment, while the other phenotype
is favoured by a glucose rich environment. This was noted in the experiments reported by Novick and Weiner, as well as others.

Second issue relates to the very essence of hysteresis --- {\em
rate-independence}. The rate-independence of the state-input
relationship in non-ideal relay and other hysteresis operators
means that the state does not depend on the rate at which the
input may have varied, but rather on the past values of the input
extrema \cite{alexei}.
Although rate-independence is an idealization, it is a useful one
because it describes an important form of memory that can not be
attained by linear dynamical systems whose memory is typically
associated with certain characteristic times, rather than input
features such as extrema.

The third point relates specifically to modeling switching of
phenotype in bacteria using the non-ideal relay operator and, in
particular, to the role and values of switching thresholds $x$ and
$y$, which quantify the response of bacteria to varying
environmental conditions and determine the bi-stability range. In
the hybrid linear differential model of M. Thattai and A. van
Oudenaarden \cite{bac1}, the maximal fitness (measured by the net
population growth rate) was achieved by the responsive switching
strategy, whereby all the bacteria switch to the currently most
favoured phenotype, provided that the rate of transitions between
phenotypes is higher than the rate of environment variations. This
strategy is described by the memoryless relay with coinciding
thresholds $x=y$ and no bi-stability region (a shifted Heaviside
step function). When the rate of transitions between phenotypes is
comparable to, or lower than, the rate of variations of the
environment, the maximal fitness in this model can be achieved by
a heterogeneous population implementing another memoryless
strategy, where bacteria anticipate fluctuations of the
environment by having a subpopulation ready in an appropriate
phenotype before the environment changes to a state favouring this
phenotype. The optimal heterogeneous distribution between
phenotypes is obtained dynamically by allowing some positive rate
of transitions from the currently most favoured phenotype to an
unfavoured one. E. Kussell and S. Lieber argued that there is a
cost of maintaining sensory machinery to respond to changes of the
environment and introduced a cost of sensing  in the form of an
explicit reduction in the growth rate  \cite{bac2}. In their
linear differential model, random switching of phenotype, which is
not correlated to slowly varying environmental conditions, can
confer more fitness to the population than the responsive
switching when the penalty for sensing is large.  A modification
of these models, which incorporates a natural switching cost in
the form of a temporary inhibition of the reproductive activity in
bacteria undergoing a transition to a different phenotype (a lag
phase), shows that hysteresis characterized by two different
switching thresholds $x>y$ can also confer fitness to bacteria
\cite{we}\footnote{{\em In vitro} experimental studies give
evidence that the process of changing phenotype is stressful for
bacteria. In particular, bacteria may not reproduce within a
period of time preceding, during, or following this process. In
these experiments, a colony of bacteria grown in a Petri dish with
one nutrient is swapped to a Petri dish with another nutrient.
After a period of inactivity, or a shock, following the swap,
bacteria start a transition to the other phenotype which is better
fit for consuming the new type of food.}. In the {\em adiabatic
limit} of slowly varying environmental conditions, the switching
pattern of bacteria in the model is described by the hysteretic
non-ideal relay shown in Figure \ref{figure1}. In faster uncertain
environments, the maximal growth was shown to be the one that
results from a form of stochastic resonance where  the internal
characteristic time associated with the growth delay is about
equal to the characteristic time between subsequent phenotype
switching events (the latter time is defined by the length of the
bi-stability interval $(y,x)$ and the characteristic time of the
environment variations). More frequent switching in faster
environments causes the organisms to keep delaying their
reproduction; slower environments causes the bacteria to wait too
long in a sub-optimally reproducing phenotype.

An advantage of using strategies with memory has been also shown in a different framework, namely, the
game theory \cite{g1} applied to models where bacteria are considered as players
in an evolutionary game \cite{bac2, h2, c12}.

The idea of diversification or bet-hedging has been discussed in different biological contexts in many publications, often without reference to any specific mechanism by which it can be implemented. For example, the view that diversity (heterogeneity) can help
improve fitness in varying conditions is very well established in ecology. The main idea behind diversification in epigenetics is that genetically identical organisms can  grow their numbers faster by judiciously choosing a certain proportion of their population to be in a currently less favored phenotype when the environmental varies in time. In particular, it has been shown that most appropriate random phenotype choice strategies are based on selecting phenotype switching probabilities that are tuned in some sense to the environmental uncertainties. The main finding of \cite{bac1} was that bet-hedging in
the form of the anticipating switching strategy described above can lead to the maximal growth rate. If a population of bi-stable organisms includes
subpopulations with different bi-stability ranges $(y,x)$, then bet-hedging can be realized by allowing permanently ongoing transitions between
subpopulations. Due to the difference in the bi-stability range, each subpopulation is tuned to a specific pattern of variation of the environment, which maximizes
the growth rate of this subpopulation. An exchange process between the subpopulations can bet-hedge against changes of this pattern, at least in principle.

\section{Model description}\label{secModelDescription}

\subsection{Objective and modeling assumptions}
In this paper, we attempt to formulate a class of models, which
account for a number of phenomena listed above, namely (a)
switching of bacteria between two phenotypes in response to
variations of environmental conditions; (b) hysteretic switching
strategy (switching rules) associated with bi-stability of
phenotype states; (c) heterogeneity of the population in the form
of a distribution of switching thresholds; (d) bet-hedging in the
form of diffusion between subpopulations characterized by
different bi-stability ranges; and, (e) competition for nutrients.
The resulting model is a reaction-diffusion system including, as
reaction terms, discontinuous hysteresis relay operators and the
integral of those. This integral can be interpreted as the
Preisach operator \cite{alexei} with a time dependent density (the
density is a component of the solution describing the varying
distribution of bacteria). The main objective of this paper is to
prove well-posedness of the model. In the last part of the paper
we will also present a preliminary discussion of a few dynamic
scenarios obtained numerically, where fitness, competition and
diffusion act together to select a certain distribution of
switching thresholds in the population. In particular, the model
seems to demonstrate interesting dynamics such as pattern
formation and oscillations. However, more detailed and rigorous
analysis of dynamics is beyond the scope of this paper and remains
the subject of future work.

We assume that each of the two phenotypes, denoted by $1$ and $-1$, consumes a different type of nutrient
(for example, one consumes lactose, the other consumes glucose).
The amount of nutrient available for phenotype $i$ at the moment $t$ is denoted by $f_i(t)$ where $i=\pm 1$.
We base our model on the following assumptions.

\begin{itemize}
\item Each bacterium changes phenotype in response to the variations of the variable $w=f_1/(f_1+f_{-1})-1/2$.
\end{itemize}
This  variable measures the deviation of the relative concentration
of the first nutrient from the value 1/2 in the mixture of the two nutrients.
Bacteria sense changes of the nutrient concentrations and, by changing
to the phenotype for which more food is available, can potentially increase the growth rate
of the population.

\begin{itemize}
\item The input $w=w(t)$ is mapped to the binary time trace
$r=r(t)$ of the phenotype (state) of a bacterium by the non-ideal
relay operator ${\mathcal R^x}$ (see Figure \ref{figure1}) with
symmetric switching thresholds $x$ and $y=-x$, where $x>0$. This
binary function of time will be denoted by $r={\mathcal R^x}(w)$.
\end{itemize}
As transitions between the states of a relay are instantaneous,
this assumption relates to the situation where the rate of
transitions between phenotypes is much higher than the rate of the
input variations. The non-ideal relay operator ${\mathcal R^x}$
is rigorously defined in the next subsection.

\begin{itemize}
\item The population includes bacteria with different bi-stability
ranges $(-x,x)$; the threshold value $x$ varies over an interval
$[\underline{x}, \overline{x}]$, where $0<\underline{x} <
\overline{x}$.
\end{itemize}

We will denote by $u(x,t)$ the density of the biomass of bacteria with given switching thresholds $\pm x$ at a moment $t$.
\begin{itemize}
\item There is a diffusion process acting on the density $u$.
\end{itemize}
The diffusion process models sporadic changes of the switching threshold $x$ in bacteria. Another source of diffusion can be the birth process if we assume that a bacterium with a threshold $x$ produces offsprings with different thresholds, for example, according to the Gaussian distribution centered at $x$.
The diffusion can be viewed as a bet-hedging mechanism in the population.

Finally, we make the following important assumption.
\begin{itemize}
\item
   At any particular time moment $t$, for any given $x$, all the bacteria with the switching threshold values $\pm x$ are in the same state (phenotype).
\end{itemize}
   That is, $u(x,t)$ is the total density of bacteria with the threshold $x$ at the moment $t$ and they are all in the same state.
It means that when a bacterium with a threshold $x'$
sporadically changes its threshold to a different value $x$, it simultaneously
copies the state from other bacteria which have the threshold $x$. (Or, if a bacterium with a threshold $x'$ produces
an offspring with a threshold $x$, the offspring copies the state of other bacteria with the threshold $x$). In particular,
this may require a bacterium to change the state when its threshold changes.
Models where the state of a bacterium remains unchanged after a change of the threshold will be considered in a different work.

With these assumptions, we obtain the following model of the
evolution of bacteria and nutrients:
\begin{equation}\label{00}
\left\{
\begin{aligned}
& {u}_t=  u_{xx} + \frac12  (1+{\mathcal R}^x(w))\, u
f_1+\frac12(1-{\mathcal R}^x(w))\, u f_{-1},\\
& \dot{f}_{1} =-\frac12
f_1\int_{\underline{x}}^{\overline{x}}(1+{\mathcal R}^x(w))\,
u\,dx,\\
& \dot{f}_{-1} =-\frac12
f_{-1}\int_{\underline{x}}^{\overline{x}}(1-{\mathcal R}^x(w))\,
u\,dx,
\end{aligned}
\right.
\end{equation}
where $u_t$ and $u_{xx}$ are the derivatives of the population
density $u$, dot denotes the derivative with respect to time, and
all the non-ideal relays ${\mathcal R^x}$,
$x\in[\underline{x},\overline{x}]$, have the same input
$w=f_1/(f_1+f_{-1})-1/2$.
 Here we additionally assume the growth rate $\frac12 (1+ i {\mathcal R}^x(w))\, u f_i $ based on the mass action law for bacteria in the phenotype $i=\pm 1$. This growth rate is proportional to the product of the population density $u$ and the concentration $f_i$ of the nutrient preferred by the phenotype $i$ with the coefficient of proportionality scaled to unity. The diffusion coefficient is set to unity. The rate of the consumption of nutrient in the equation for $f_i=f_i(t)$ is proportional to the total biomass of bacteria in the phenotype $i$
(i.e., the biomass of all bacteria eating this type of nutrient),
hence the integral  (the coefficient of proportionality is also
set to unity for simplicity); $\underline{x}$ and $\overline{x}$
are the lower and upper bounds on available threshold values,
respectively.

We assume that a certain amount of nutrients is available at the
initial moment; the nutrients are not supplied after that moment.
 Bacteria do not die but stop growing when all the nutrient has been consumed.
We assume the Neumann boundary conditions for $u$, i.e., no flux
of the population density $u$ through the lower and upper bounds
of available threshold values.


\subsection{Rigorous model setting}\label{secSetting}

Throughout the paper,  we assume  that $x\in[\ux,\ox]$, where
$0<\ux<\ox$.


We begin with a rigorous definition of the hysteresis operator
(non-ideal relay) with fixed thresholds $\pm x$. The {\em
non-ideal relay} is the  operator $\cR^x$ which takes continuous
functions $w=w(t)$ defined on an interval $[0,T)$ to the binary
functions $r=\cR^x(w)$ of time defined on the same interval. The
operator $\cR^x$ is given by
\begin{equation}\label{nonideal}
 {\mathcal R}^x (w)(t) = \left\{
\begin{array}{rll} -1 & {\rm if} & w(\tau)< -x \ {\rm for\ some} \ \tau\in [0,t] \\
 && {\rm and} \ w(s)<x \ {\rm for\ all} \ s\in [\tau,t],
 \\1& {\rm if} &  w(\tau)\ge x \ {\rm for\ some} \ \tau\in [0,t] \\
 && {\rm and} \ w(s)>-x \ {\rm for\ all} \ s\in [\tau,t],\\
 r_0 & {\rm if} & -x \le w(\tau) < x \ {\rm for\ all} \ \tau\in
 [0,t],
 \end{array}
\right.
\end{equation}
where $r_0$ is either $1$ or $-1$ ({\em initial configuration} of
the non-ideal relay $\cR^x$). Since $r_0$ may take different
values for different $x$, we write $r_0=r_0(x)$, where $r_0(x)$ is
a given function taking values $\pm 1$. In what follows, we  do
not explicitly indicate the dependence of the operator ${\mathcal
R}^x$ on $r_0(x)$. Some additional assumptions on the structure of
the initial configuration $r_0(x)$ will be made in
Sec.~\ref{secStatesPreisach} (in particular, see
relations~\eqref{eqCompatibility}).

The {\em distributed relay operator} $\cR(w)$ taking functions
$w=w(t)$ to functions $r=r(x,t)$ is defined by
$\cR(w)(x,t)={\mathcal R}^x(w)(t)$.

We also denote
\begin{equation}\label{eqNP}
\cU(u)(t)=\int\limits_{\underline{ x}}^{\overline{ x}} u( x,t)\,d
x,\quad \cP(u,w)(t)=\int\limits_{\underline{ x}}^{\overline{ x}}
u( x,t) \cR^x(w)(t)\,d x
\end{equation}
and call $\cP$ the {\em Preisach operator} (see more details in
Sec.~\ref{secStatesPreisach}).

Now let us replace the unknown function $f_1$ and $f_{-1}$ in
system \eqref{00} with $v=f_1+f_{-1}$ and
$w=f_1/(f_1+f_{-1})-1/2$. The resulting system has the form
\begin{equation}\label{eqBacteriaGeneral}\left\{
\begin{aligned}
& u_t = u_{xx} +\cA(v,w)u,\\
&\dot v=\cV(u,v,w),\\
&\dot w=\cW(u,w),
\end{aligned}\right.
\end{equation}
where we assume the Neumann
boundary conditions
\begin{equation}
\label{eqBC}  u_x|_{ x=\ux} = u_x|_{ x=\ox}=0
\end{equation}
and the initial conditions
\begin{equation}
\label{eqIC} u( x,0)=u_0( x),\quad v(0)=v_0,\quad w(0)=w_0.
\end{equation}
Here
$$
\begin{aligned}
&\cA(v,w)=\left(\frac12 + w
 \cR(w)\right)v,\\
&\cV(u,v,w)=-\left(\frac{1}{2}  \cU(u)+w
\cP(u,w)\right)v,\\
&\cW(u,w)=- \left(\dfrac{1}{2}+w\right)\left(\dfrac{1}{2}-w\right)
\cP(u,w),
\end{aligned}
$$
and the initial configuration $r_0(x)$ of the distributed relay
operator $\cR$ is a part of initial data.


Let $T>0$, and let  $Q_T=(\ux,\ox)\times(0,T)$. We will use the
following spaces:
\begin{enumerate}
\item The standard Lebesgue spaces $L_ q(Q_T)$ and $L_ q=L_
q(\ux,\ox)$ with $1\le q\le\infty$.

\item The Sobolev spaces $W_ q^k=W_ q^k(\ux,\ox)$, $k\in\bbN$.

\item The anisotropic Sobolev space $W_2^{2,1}(Q_T)$  with the
norm
$$
\|u\|_{W_2^{2,1}(Q_T)}=\left(\int\limits_0^T
\|u(\cdot,t)\|_{W_2^2}^2\, dt+ \int\limits_0^T
\|u_t(\cdot,t)\|_{L_2}^2\, dt\right)^{1/2}.
$$

\item The space $C([0,T];L_ q)$ of $L_ q$-valued functions
continuous in $t\in[0,T]$.
\end{enumerate}

\section{Preisach model}\label{secStatesPreisach}

\subsection{States of the Preisach model}

In this section, we establish some continuity properties of the
operators $\cR(w)$ and $\cP(u,w)$ (the latter defined
by~\eqref{eqNP}). In particular, we will show that, for any $w\in
C[0,T]$ with bounded variation and $u\in (C[0,T];L_1)\cap
L_\infty(Q_T)$, the function $\cR(w)$ belongs to $C([0,T];L_ q)$
with any $ q\ge 1$ and the function $\cP(u,w)$ belongs to
$C[0,T]$.

We note that the operator $\cP(u,w)$  is a specific case of the {\em Presiach
operator with density $u(x,t)$ and input $w(t)$} \cite{alexei}. However, in the
literature, the Preisach operator is usually considered for time independent density.

We begin with some definitions. Fix $T>0$.

\begin{definition} For a given input $w\in C[0,T]$, the {\em state} of the Preisach model at the moment $t\in [0,T]$ is defined as a
subset $A(t)$ ($t\in[0,T]$) of the interval $[\underline
x,\overline x]$ given by
$$
A(t)=\{x\in [\underline x,\overline x] : \cR^x(w)(t)=1\}.
$$
\end{definition}

\begin{definition} For a given input $w\in C[0,T]$, we say that the state $A(t)$ is {\em simple} for some $t\in[0,T]$
if it is a union of finitely many disjoint intervals $[ x_k^-,
x_k^+]$, $k=1,\dots, N$.
\end{definition}

In what follows, we assume that $A(0)$ is simple and satisfies the
compatibility condition
\begin{equation}\label{eqCompatibility}
[\underline  x, w(0)]\subset A(0) \quad \text{if }   w(0)\ge
\underline x; \qquad [\underline  x,-w(0))\cap
A(0)=\varnothing\quad \text{if }   w(0)< -\underline x.
\end{equation}
We will show that $A(t)$ remains simple and satisfies the
compatibility condition
\begin{equation}\label{eqCompatibilityt}
[\underline  x, w(t)]\subset A(t) \quad \text{if }   w(t)\ge
\underline x; \qquad [\underline  x,-w(t))\cap
A(t)=\varnothing\quad \text{if }   w(t)< -\underline x
\end{equation}
for all $t\in[0,T]$.

We remind that we have agreed that each relay $R_x(w)$ switches to
the state $1$ at the moment when the input reaches the value $ x$
and remains in the same state $1$ when the input reaches the value
$- x$, switching back to the state $-1$ only when the input
becomes smaller than $- x$. However, this is not essential as
including/not including the end points of $[ x_k^-, x_k^+]$ in
$A(t)$ is not important (one cares about the classes of sets $A$
which coincide almost everywhere).


Suppose, at some moment $t_1\in[0,T)$, the state is $A(t_1)$. Fix
some moment $t_2 \in(t_1,T]$ and consider $A(t)$ for
$t\in[t_1,t_2]$. From the definition of dynamics of individual
relays, it follows that the state $A(t)$ will evolve in response
to the changes of a continuous input $w(t)$ after the moment $t_1$
according to the following rules.

Consider the running maximum and minimum of $w$:
$$
M(t,t_1)=\max_{s\in [t_1,t]} w(s),\qquad m(t,t_1)=\min_{s\in
[t_1,t]} w(s).
$$

\begin{enumerate}

\item If the input satisfies $-\underline x \le w(t) $ on the time
interval $[t_1, t_2]$, then on this time interval
\begin{equation}\label{bu}
A(t)= A(t_1) \cup J(t),
\end{equation}
where
$$
J(t)=\varnothing \quad {\rm if} \quad M(t,t_1)<\underline x;
\qquad J(t)=[\underline x, M(t,t_1)] \quad {\rm if} \quad
M(t,t_1)\ge \underline x.
$$
In particular, $A(t)$ is expanding (non-strictly). We also note
that, during any time interval within which $-\underline x\le
w(t)<\underline  x$, the set $A(t)$ does not change.

\item If the input satisfies $w(t)<\underline x$ on the time
interval $[t_1, t_2]$, then on this time interval
\begin{equation}\label{bv}
A(t)= A(t_1) \setminus I(t)
\end{equation}
where
$$
I(t)=\varnothing \quad \text{if } m(t,t_1)\ge -\underline x; \quad
I(t)=(\underline x,-m(t,t_1)) \quad \text{if }  m(t,t_1)<
-\underline x.
$$
In this case, $A(t)$ is contracting (non-strictly). Again, $A(t)$
does not change during any time interval within which $-\underline
x\le w(t)<\underline  x$.

\end{enumerate}

Now we consider the evolution of the state $A(t)$ for all
$t\in[0,T]$, provided that the input $w(t)$ is continuous and has
a bounded variation. Specifically, we assume that
$$
{\rm Var}_{0}^T [w]\le K
$$
for some $K=K(T)>0$. Such inputs form a closed convex subset of
$C[0,T]$ for any given $T>0$. For every such input, there is a
finite sequence of moments $0<t_1<t_2<\cdots<T$ such that, on any
interval $[t_k,t_{k+1}]$, either $-\underline x \le m(t,t_k) $, or
$M(t,t_k)<\underline x$, or both these relations hold
simultaneously. Hence, on each $[t_k,t_{k+1}]$, at least one of
the above two rules of evolution of the state $A(t)$ applies, thus
defining the state $A(t)$ at any moment  $t\in[0,T]$ by a finite
number of applications of formulas (\ref{bu}), (\ref{bv}) on
successive time intervals $[t_k,t_{k+1}]$. We note that the
sequence $t_k$ is not unique; we call such sequences {\em
admissible partition} sequences. The     dynamics of $A(t)$,
according to the above rules, is independent of the choice of an
admissible partition sequence.

Note that, according to the definition of the relay, the states
$A(t)$ are simple and satisfy the compatibility
condition~\eqref{eqCompatibilityt} for any $t\in[0,T]$.

\subsection{Continuity}\label{secLipRP} For any
measurable sets $B_1,B_2\subset[\underline x,\overline x]$, denote
$$
\rho(B_1,B_2)={\rm meas}\,\{(B_1\setminus B_2)\cup (B_2\setminus
B_1)\}.
$$

Consider the evolution of the states $A_1(t)$ and $A_2(t)$ in
response to the inputs $w_1(t)$ and $w_2(t)$, respectively, for
$t\in[0,T]$.

\begin{lemma}\label{lLipContRho}
Let $w_1,w_2\in C[0,T]$ for some $T>0$, and let
$$
{\rm Var}_{0}^T [w_j]\le K,\quad  j=1,2,
$$
with some $K=K(T)$. Let $A_j(t)$ be the state corresponding to the
input $w_j(t)$. If $A_j(0)$, $j=1,2$, are simple and satisfy the
compatibility condition~\eqref{eqCompatibility}, then
\begin{equation}\label{lipsc}
\max_{t\in[0,T]}\rho(A_1(t),A_2(t))\le \rho(A_1(0),A_2(0))+ L
\|w_1-w_2\|_{C[0,T]},
\end{equation}
where $L = 2+K /(2\underline x)$.
\end{lemma}
\proof First, we choose a finite sequence $0<t_1<t_2<\dots<T$
which is an admissible partition sequence for the evolution of
both states $A_1(t)$ and $A_2(t)$ simultaneously.

To do so, we fix an arbitrary number $x_0\in(0,\underline x)$ and
consider the auxiliary relay $\cR^{x_0}$.
 Define the   sequence of all moments $ t_1<t_2<\dots $ on the interval $[0,T]$ as
the moments when the output  $\cR^{x_0}(w_1)(t)$ switches between
the states $\pm1$ in response to the input $w_1(t)$. This sequence
is finite since $w_1$ is of bounded variation. Then
$0,t_1,t_2,\dots,T$ is an admissible partition sequence for the
evolution of the state $A_1(t)$. Assume that, on the time interval
$[0, T]$,
$$
\|w_1-w_2\|_{C[0,T]} \le \underline  x-x_0.
$$
Then, the same sequence $0,t_1,t_2,\dots,T$ is an admissible
partition sequence for the evolution of the state $A_2(t)$ on the
interval $[0,T]$. Moreover, the first updating rule applies to
both evolutions on each time interval $[t_k,t_{k+1}]$ where
$\cR^{x_0}(w_1)(t)=1$, while the second updating rule applies to
both evolutions on each time interval $[t_k,t_{k+1}]$  where
$\cR^{x_0}(w_1)(t)=-1$. Applying formulas~\eqref{bu}
and~\eqref{bv} to both evolutions, we see that
$$
\rho(A_1(t),A_2(t))\le \rho(A_1(t_k),A_2(t_k))+
\|w_1-w_2\|_{C[0,T]},\qquad t\in [t_k,t_{k+1}].
$$
 As the number of switching points $t_k$ in the
interval $[0,T]$ does not exceed the value $1+{\rm Var}_{0}^{T}
[w_1]/(2 x_0)\le 1+K/(2x_0)$, we obtain the Lipschitz estimate
$$
\max_{t\in[0,T]}\rho(A_1(t),A_2(t))\le \rho(A_1(0),A_2(0))+ L
\|w_1-w_2\|_{C[0,T]}
$$
with $L =2+K/(2 x_0)$. Since $x_0\in(0,\underline x)$ is
arbitrary, we have~\eqref{lipsc}.
\endproof

Using Lemma~\ref{lLipContRho}, we obtain the main results of this
section about the continuity of the operators $\cR(w)$ and
$\cP(u,w)$ (the latter defined by~\eqref{eqNP}.

\begin{lemma}\label{lContRP}
\begin{enumerate}
\item Let $w\in C[0,{T}]$ and $\Var_0^{T}[w]$ be finite. Then
$\cR(w)\in C([0,{T}];L_ q)$ with any $ q\ge 1$.

\item If, additionally, $u\in C([0,{T}];L_1)\cap L_\infty(Q_{T})$,
then $\cP(u,w)\in C[0,{T}]$.
\end{enumerate}
\end{lemma}

\proof 1. Denote $r(x,t)=\cR^x(w)(t)$. Then, for any fixed
$t_0\in[0,T]$,
$$
\int\limits_{\ux}^{\ox}|r(x,t)-r(x,t_0)|^qdx=2^q
\rho(A(t),A(t_0))\to 0,\quad t\to t_0,
$$
due to the updating rules~\eqref{bu} and~\eqref{bv}. This proves
assertion 1.

2. To prove the continuity of the function $\cP(u,w)(t)$, we
estimate
$$
\begin{aligned}
&\int\limits_{\ux}^{\ox}|u(x,t)r(x,t)-u(x,t_0)r(x,t_0)|dx \\
\quad &\le
\int\limits_{\ux}^{\ox}|u(x,t)-u(x,t_0)|\,|r(x,t)|dx+\int\limits_{\ux}^{\ox}|u(x,t_0)|\,|r(x,t)-r(x,t_0)|dx\\
\quad &\le
\|u(\cdot,t)-u(\cdot,t_0)\|_{L_1}+2\|u(\cdot,t_0)\|_{L_\infty}\,\rho(A(t),A(t_0))\to
0,\quad t\to t_0,
\end{aligned}
$$
due to the assumptions of the lemma and the updating
rules~\eqref{bu} and~\eqref{bv}.
\endproof

\begin{lemma}\label{lLipRP}
\begin{enumerate}
\item Let
$$
w_j\in C[0,{T}],\quad \Var_0^{T}[w_j]\le K,\quad A_1(0)=A_2(0),
$$
where $K>0$ $(j=1,2)$.  Then, for any $ q\ge 1$,
$$
\|\cR(w_1)-\cR_2(w)\|_{C([0,{T}];L_ q)} \le
L_\cR\|w_1-w_2\|_{C[0,{T}]}^{1/ q},
$$
where $L_\cR=L_\cR(K, q)>0$.

\item If, additionally,
$$
u_j\in C([0,{T}];L_1),\quad \|u_j\|_{L_\infty(Q_{T})}\le c
$$
for some $c>0$ $(j=1,2)$, then
$$
\|\cP(u_1,w_1)-\cP(u_2,w_2)\|_{C([0,{T}])}\le
L_\cP(\|u_1-u_2\|_{C([0,{T}];L_1)}+\|w_1-w_2\|_{C[0,{T}]}),
$$
where $L_\cP=L_\cP(c,K)>0$.
\end{enumerate}
\end{lemma}
\proof 1. Denote $r_j(x,t)=\cR^x(w_j)(t)$, $j=1,2$. Then, using
Lemma~\ref{lLipContRho}, we have
$$
\int\limits_{\ux}^{\ox}|r_1(x,t)-r_2(x,t)|^qdx=2^q
\rho(A_1(t),A_2(t))\le 2^q L\|w_1-w_2\|_{C[0,T]}.
$$

2. To prove the continuity of the operator $\cP$, we estimate
$|\cP(u_1,w_1)(t)-\cP(u_2,w_2)(t)|$ as follows (omitting the
arguments of the integrands):
$$
\begin{aligned}
 \int\limits_{\ux}^{\ox}|u_1 r_1-u_2 r_2|dx  &\le
\int\limits_{\ux}^{\ox}|u_1|\,|r_1-r_2|dx+\int\limits_{\ux}^{\ox}|u_1-u_2|\,|r_2|dx\\
 &\le 2 c \rho(A_1(t),A_2(t))+
\|u_1(\cdot,t)-u_2(\cdot,t)\|_{L_1}\\
&\le 2cL \|w_1-w_2\|_{C[0,T]} + \|u_1-u_2\|_{C([0,T];L_1)},
\end{aligned}
$$
where Lemma~\ref{lLipContRho} was used to estimate
$\rho(A_1(t),A_2(t))$.
\endproof

\begin{remark}\label{remNonlipR}
We underline (see the proof of Lemma~\ref{lLipRP}) that the
operators
$$
\cR^x: C[0,T]\to L_q,\quad \cR: C[0,T]\to C([0,T];L_q)
$$
are continuous for any $q\ge 1$. However, they are Lipschitz
continuous {\em only} for $q=1$.
\end{remark}

\section{Well-posedness}\label{secWellPosedness}

In this section, we establish existence (first locally and then
globally) and uniqueness  for
problem~\eqref{eqBacteriaGeneral}--\eqref{eqIC}.

We will often write  $U$ or $U(t)$ instead of $\cU(u)(t)$
(see~\eqref{eqNP}). This should lead to no confusion. We will also
denote
$$
U_0=\cU(u)(0)=\int\limits_{\ux}^{\ox} u_0( x)\,d x.
$$

 Since
the right-hand sides in~\eqref{eqBacteriaGeneral} contain the
relays $\cR^x$, one has to fix the initial state of those relays.
We assume throughout that the initial state $A(0)$ is simple and
satisfies the compatibility condition~\eqref{eqCompatibility}.

\subsection{Linear parabolic problem}

In this subsection, we formulate some auxiliary results on the
following linear parabolic problem:
\begin{equation}\label{eqLinearParabolic}
\left\{
\begin{aligned}
& u_t = u_{xx} + a( x,t)u+f( x,t),\quad  x\in(\ux,\ox),\ t>0,\\
& u( x,0)=u_0( x)
\end{aligned}\right.
\end{equation}
with the homogeneous Neumann boundary conditions. In what follows,
we will use the functional spaces defined in
Sec.~\ref{secSetting}.

The first lemma follows from~\cite[Chap. 4]{LadSolUral}.
\begin{lemma}\label{lLinearParabolic1}
Let $T\le T_0$ for some $T_0>0$, and let $\|a\|_{L_\infty(Q_T)}\le
\oa$ for some $\oa\ge 0$. Let $f\in L_2(Q_T)$ and $u_0\in W_2^1$.
Then problem~\eqref{eqLinearParabolic} has a unique solution $u\in
W_2^{2,1}(Q_T)$ and
$$
\|u\|_{W_2^{2,1}(Q_T)}\le c(\|u_0\|_{W_2^1}+\|f\|_{L_2(Q_T)}),
$$
where $c=c(\oa,T_0)>0$ does not depend on $u$, $a(x,t)$ and $T$.
\end{lemma}

The second lemma deals with   continuous dependence of solutions
on the coefficient $a(x,t)$. Consider a sequence $a_j\in
L_\infty(Q_T)$, $j=1,2,\dots$. Denote by $u_j$ the solution of
problem~\eqref{eqLinearParabolic} with $a_j$ instead of $a$.

\begin{lemma}\label{lLinearParabolic2}
Let  $\|a_j\|_{L_\infty(Q_T)}\le \oa$ and $\|a-a_j\|_{L_2(Q_T)}\to
0$ as $j\to\infty$. Then
$$
\|u-u_j\|_{W_2^{2,1}(Q_T)}\to0,\quad j\to\infty.
$$
\end{lemma}
\proof The function $m_j=u-u_j$ is a solution of the problem
$$
\left\{
\begin{aligned}
& m_{jt} = m_{jxx} + a_j( x,t)m_j+(a-a_j)u,\\
& u( x,0)=0.
\end{aligned}\right.
$$
Therefore, by Lemma~\ref{lLinearParabolic1} and by the boundedness
of the embedding $W_2^{2,1}(Q_T)\subset L_\infty(Q_T)$, we have
$$
\|m_j\|_{W_2^{2,1}(Q_T)}\le k_1\|(a-a_j)u\|_{L_2(Q_T)}\le k_2
\|(a-a_j)\|_{L_2(Q_T)}\|u\|_{W_2^{2,1}(Q_T)},
$$
where $k_1,k_2>0$ depend  only on $\oa$ and $T$ and do not depend
on $j$. Hence, $\|m_j\|_{W_2^{2,1}(Q_T)}\to0$ as $j\to\infty$.
\endproof

\subsection{Local existence of solutions}

We introduce the space
$$
\cW(Q_T)=W_2^{2,1}(Q_T)\times C^1[0,T]\times C^1[0,T].
$$

\begin{definition}
We say that $\bu=(u,v,w)\in\cW(Q_T)$ is a ({\em strong}) {\em
solution} of problem~\eqref{eqBacteriaGeneral}--\eqref{eqIC} $($on
the interval $(0,T)$$)$ with initial data $\bu_0=(u_0,v_0,w_0)\in
W_2^1\times\bbR^2$ if
\begin{enumerate}
\item
 $u$ satisfies the first equation
in~\eqref{eqBacteriaGeneral} a.e. in $Q_T$ and the boundary
conditions~\eqref{eqBC} and the first initial conditions
in~\eqref{eqIC} in  the sense of traces.

\item $v$ and $w$ satisfy the second and the third equations
in~\eqref{eqBacteriaGeneral} and the second and the third initial
conditions in~\eqref{eqIC}, respectively, in the classical sense.
\end{enumerate}
\end{definition}

In what follows, we will often say ``solution'', meaning ``strong
solution''.

In this subsection, we will prove the following result on the
local existence of solutions.

\begin{theorem}\label{thLocalExist}
Let $\bu_0=(u_0,v_0,w_0)\in W_2^1\times\bbR^2$ and
$$
\|u_0\|_{W_2^1}\le\ou,\quad |v_0|\le \ov,\quad |w_0|\le \ow
$$
for some $\ou,\ov,\ow>0$. Then there is ${t_0}\in(0,1]$ such that
problem~$\eqref{eqBacteriaGeneral}$--$\eqref{eqIC}$ has a solution
$\bu=(u,v,w)\in\cW(Q_{t_0})$. The number ${t_0}$ depends on
$\ou,\ov,\ow$ but does not depend on $u_0,v_0,w_0$.
\end{theorem}

The idea of the proof is to construct a mapping $(v,w)\mapsto
(\tilde v,\tilde w)$ as follows.
\begin{enumerate}
\item[Step 1.]
 Given the functions $v$ and $w$, we solve the linear
parabolic problem
\begin{equation}\label{eqLinearParabolicn}
\left\{
\begin{aligned}
& u_t = u_{xx} + \cA(v,w)  u,\quad  x\in(\ux,\ox),\ t>0,\\
& u( x,0)=u_0( x)
\end{aligned}\right.
\end{equation}
with the homogeneous Neumann boundary conditions.

\item[Step 2.] After finding $u$, we find $(\tilde v,\tilde w)$ by
solving the equations
\begin{equation}\label{eqLinearEq}
\left\{
\begin{aligned}
&\frac{d \tilde v}{d t}=\cV(u,v,w),\quad \tilde v(0)=v_0,\\
&\frac{d \tilde w}{d t}=\cW(u,w),\quad \tilde w(0)=w_0.
\end{aligned}\right.
\end{equation}
\end{enumerate}
Then, using the Schauder fixed point theorem, we show that the
mapping $(v,w)\mapsto (\tilde v,\tilde w)$ has a fixed point,
which yields the solution of the original
problem~$\eqref{eqBacteriaGeneral}$--$\eqref{eqIC}$.

\begin{remark}
We note that the parabolic problem~\eqref{eqLinearParabolicn} is
well posed in $L_q$-spaces with $q>1$, while the hysteresis
operator $\cR$ is not Lipschitz continuous in these spaces  (see
Remark~\ref{remNonlipR}). Therefore, the constructed mapping
$(v,w)\mapsto (\tilde v,\tilde w)$ will be continuous, but not
Lipschitz continuous. This is the reason why we apply the Schauder
fixed point theorem and not the contraction mapping principle for
the proof of Theorem~\ref{thLocalExist}. As a result, only the
existence of solutions is proved. The uniqueness will be proved in
Sec.~\ref{secUniqueness} by using semigroups and additional
estimates in $L_1$.
\end{remark}

Let us formalize the above scheme.

We introduce the set
$$
\begin{aligned}
 B[0,{t_0}]=\{&(v,w)\in C[0,{t_0}]\times C[0,{t_0}]:\\
&\|v\|_{C[0,{t_0}]}\le 2\ov,\ \|w\|_{C[0,{t_0}]}\le 2\ow,\
\Var_0^{t_0} [w]\le Vt_0\},
\end{aligned}
$$
where $\ov$ and $\ow$ are the constants from the assumption in
Theorem~\ref{thLocalExist}, while $V>0$ and ${t_0}\le 1$ will be
chosen later on. For $(v,w)\in B[0,{t_0}]$, we will denote
$$
\|(v,w)\|_{B[0,{t_0}]}=\|v\|_{C[0,{t_0}]}+\|w\|_{C[0,{t_0}]}.
$$
Note that $B[0,{t_0}]$ is a closed convex set in $C[0,{t_0}]\times
C[0,{t_0}]$.

\begin{lemma}\label{lbaCont}
\begin{enumerate}
\item $\|\cA(v,w)\|_{L_\infty(Q_{t_0})} \le \oa$ for any $(v,w)\in
B[0,{t_0}]$, where
$$
\oa=2 |\ov|\left(\dfrac{1}{2}+2|\ow|\right).
$$
\item $\|\cA(v,w)-\cA(v_j,w_j)\|_{L_2(Q_{t_0})}\to 0$ whenever
$\|(v,w)-(v_j,w_j)\|_{B[0,{t_0}]}\to 0$.
\end{enumerate}
\end{lemma}
\proof Part 1 is trivial. Part 2 follows from~Lemma~\ref{lContRP}
(part 1) and Lemma~\ref{lLipRP} (part 1) with $ q=2$.
\endproof

Combining Lemmas~\ref{lLinearParabolic1}, \ref{lLinearParabolic2},
and~\ref{lbaCont}, we obtain the following result justifying Step
1 in the above scheme (recall that ${t_0}\le 1$).
\begin{lemma}\label{lLinearParabolicn}
Let $(v,w)\in B[0,{t_0}]$. Then problem~\eqref{eqLinearParabolicn}
has a unique solution $u\in W_2^{2,1}(Q_{t_0})$. Moreover,
\begin{enumerate}
\item the estimate
$$
\|u\|_{W_2^{2,1}(Q_{t_0})}\le c
$$
holds with $c=c(\ou,\ov,\ow)>0$ which does not depend on ${t_0}\le
1$ and $V>0;$ \item the mapping $ B[0,{t_0}]\ni(v,w)\mapsto u\in
W_2^{2,1}(Q_{t_0}) $ is continuous.
\end{enumerate}
\end{lemma}

Now we justify Step 2.

\begin{lemma}\label{lLinearEq}
There exist numbers ${t_0}\in(0,1]$ and $V>0$ such that the
following hold.
\begin{enumerate}
\item For any $(v,w)\in B[0,{t_0}]$ and for $u\in
W_2^{2,1}(Q_{t_0})$ defined by Lemma~$\ref{lLinearParabolicn}$,
the solution $(\tilde v,\tilde w)$ of problem~\eqref{eqLinearEq}
belongs to $B[0,{t_0}]$. Moreover,
\begin{equation}\label{eqLinearEq0}
\|\tilde v\|_{C^1[0,{t_0}]}+\|\tilde w\|_{C^1[0,{t_0}]}\le c_1,
\end{equation}
where $c_1>0$ depends on $\ou,\ov,\ow$, but does not depend on
$(v,w)\in B[0,{t_0}]$.

\item The mapping
$$
W_2^{2,1}(Q_{t_0})\times B[0,{t_0}]\ni (u,v,w)\mapsto (\tilde
v,\tilde w)\in B[0,{t_0}]
$$
is continuous.
\end{enumerate}
\end{lemma}
\proof 1. First, we note that
\begin{equation}\label{eqLinearEq1}
 \|u\|_{C([0,{t_0}];L_1)}\le k_1,\quad \|u\|_{L_\infty(Q_{t_0})}\le k_2,
\end{equation}
where $k_1,k_2>0$ depend on $\ou,\ov,\ow$, but do not depend on
${t_0}\le 1$ and $V>0$. Indeed, for any ${t_0}\le 1$, we can
extend $v(t)$ and $w(t)$ to $[0,1]$ as continuous functions
without increasing their norms and without changing the variation
of $w$. By Lemma~\ref{lLinearParabolicn}, we obtain a unique
solution $u\in W_2^{2,1}(Q_1)$ of
problem~\eqref{eqLinearParabolicn} on the time interval $(0,1)$
such that
$$
\|u\|_{W_2^{2,1}(Q_1)}\le c(\ou,\ov,\ow).
$$
This estimate and the boundedness of the embeddings
$W_2^{2,1}(Q_1)\subset C([0,1];L_1)$ and $W_2^{2,1}(Q_1)\subset
L_\infty(Q_1)$ imply~\eqref{eqLinearEq1} with $k_1,k_2>0$ not
depending on ${t_0}\le 1$.

Using Lemma~\ref{lContRP}  and estimates~\eqref{eqLinearEq1}, we
see that, for any $(v,w)\in B[0,{t_0}]$,
\begin{equation}\label{eqLinearEq2}
\|\cV(u,v,w)\|_{C[0,{t_0}]}\le k_3,\quad
\|\cW(u,w)\|_{C[0,{t_0}]}\le k_4,
\end{equation}
where $k_3,k_4>0$ depend on $\ou,\ov,\ow$, but do not depend on
${t_0}\le 1$ and $V>0$.

Now we choose
\begin{equation}\label{eqLinearEq3}
{t_0}=\min\left(\dfrac{\ov}{k_3},\dfrac{\ow}{k_4},1\right),\quad
V=k_4.
\end{equation}
Then the solution $(\tilde v,\tilde w)$ of
problem~\eqref{eqLinearEq} belongs to $B[0,{t_0}]$.

Estimate~\eqref{eqLinearEq0} follows from~\eqref{eqLinearEq2}.

2. Part 2 of the lemma follows from part 2 of Lemma~\ref{lLipRP}.
\endproof

\proof[Proof of Theorem~$\ref{thLocalExist}$] Combining
Lemmas~\ref{lLinearParabolicn} and~\ref{lLinearEq} with the
compactness of the embedding $C^1[0,t_0]\subset C[0,t_0]$ and
using the Schauder fixed-point theorem, we conclude that the
mapping $(v,w)\mapsto (\tilde v,\tilde w)$ has a fixed point,
which yields the solution of the original
problem~$\eqref{eqBacteriaGeneral}$--$\eqref{eqIC}$. \qed

\subsection{Global existence of solutions}

Our next goal is to prove that the local solution of
problem~\eqref{eqBacteriaGeneral}--\eqref{eqIC} given by
Theorem~\ref{thLocalExist} can be extended to an arbitrarily large
time interval. Here we will concentrate on the physically relevant
case where the initial data $\bu_0=(u_0,v_0,w_0)$ satisfies
\begin{equation}\label{eqPhysicalInitData}
u_0( x)\ge0,\quad v_0\ge 0, \quad |w_0|\le 1/2.
\end{equation}

First, we prove some a priori estimates  of solutions.

\begin{lemma}\label{lPositiveSolutions}
Let $\bu=(u,v,w)\in \cW(Q_T)$ be a solution of
problem~\eqref{eqBacteriaGeneral}--\eqref{eqIC} on some time
interval $(0,T)$ with initial data $\bu_0=(u_0,v_0,w_0)$
satisfying~\eqref{eqPhysicalInitData}. Then the following hold.
\begin{enumerate}
\item $u( x,t)\ge 0$ for all $( x,t)\in Q_T$.

\item $  v(t)\ge 0$ for all $t\in[0,T]$.

\item $|w(t)|\le 1/2$ for all $t\in[0,T]$.
\end{enumerate}
\end{lemma}
\proof 1. Denote $a( x,t)=\cA(v,w)$. Then $a\in L_\infty(Q_T)$,
while the first equation in~\eqref{eqBacteriaGeneral} takes the
form
\begin{equation}\label{eqPositiveSolutions1}
 u_t = u_{xx} +a( x,t)u
\end{equation}
Along with~\eqref{eqPositiveSolutions1}, we consider the equations
\begin{equation}\label{eqPositiveSolutions2}
 u^\varepsilon_t= u_{xx}^\varepsilon +a(
x,t)u^\varepsilon+\varepsilon,
\end{equation}
where $\varepsilon> 0$. Due to Lemma~\ref{lLinearParabolic1},
equation~\eqref{eqPositiveSolutions2} with the Neumann boundary
conditions and the initial condition
$$
u^\varepsilon( x,0)=u_0( x)+\varepsilon
$$
has a unique solution $u^\varepsilon\in W_2^{2,1}(Q_T)$ and
\begin{equation}\label{eqPositiveSolutions3}
\|u^\varepsilon-u\|_{W_2^{2,1}(Q_T)}\to 0,\quad \varepsilon\to 0.
\end{equation}

We fix $\varepsilon>0$. Denote
$$
b( x,t,u^\varepsilon)=a( x,t)u^\varepsilon+\varepsilon.
$$
Then we have the following:
\begin{enumerate}
\item[(a)] $b( x,t,0)=\varepsilon>0$ for $( x,t)\in \oQ_T$,

\item[(b)] $b( x,t,\cdot)$ is continuous near the origin (i.e.,
near $u^\varepsilon=0$) uniformly with respect to $( x,t)\in
\oQ_T$,

\item[(c)] $u^\varepsilon( x,0)\ge\varepsilon>0$.
\end{enumerate}
Regularizing the right-hand side $b( x,t,u_\varepsilon( x,t))$ and
applying the method of invariant regions (see~\cite{Smoller} for
classical solutions), we obtain that $u^\varepsilon( x,t)\ge 0$
for $( x,t)\in Q_T$. Hence, using \eqref{eqPositiveSolutions3} and
the continuity of the embedding $W_2^{2,1}(Q_T)\subset C(\oQ_T)$
yields $u( x,t)\ge 0$ for $( x,t)\in Q_T$.

2. The second equation in~\eqref{eqBacteriaGeneral} can be written
as
$$
\dot v=c(t)v
$$
with appropriate continuous function $c(t)$. Obviously, $v(t)\ge
0$ if $v_0\ge0$.

3. The third equation can be treated similarly to the second one.
\endproof

In the next lemma, we estimate the solutions from above. This will
allow us to prove the existence of solutions on arbitrarily large
time interval. The estimates will involve the total population
$\cU(u)(t)$ (see~\eqref{eqNP}). Whenever it appears, we keep
writing  $U$ or $U(t)$, as before.

\begin{lemma}\label{lEstimatesSolutions}
Let $\bu=(u,v,w)\in \cW(Q_T)$ be a solution of
problem~\eqref{eqBacteriaGeneral}--\eqref{eqIC} on some time
interval $(0,T)$ with initial data $\bu_0=(u_0,v_0,w_0)$
satisfying~\eqref{eqPhysicalInitData}. Then the following hold.
\begin{enumerate}
\item  $v(t)\le v_0$ for all $t\in[0,T]$.

\item   $U(t)+v(t)= U_0+v_0.$

\item There is a function $\bar u(t)$, $t\ge 0$, depending on
$\bu_0$, but not on~$T$, bounded on bounded sets and such that
$$
\|u(\cdot,t)\|_{W_2^1}\le \bar u(t).
$$
\end{enumerate}
\end{lemma}
\proof 1. By Lemma~\ref{lPositiveSolutions}, $u( x,t)\ge0$ and
$|w(t)|\le 1/2$. Therefore,
$$
\dfrac{U(t)}{2}+w(t)\cP(u,w)(t)\ge 0.
$$
Hence, the second equation
in~\eqref{eqBacteriaGeneral} yields $ \dot v\le 0, $ which implies
the first assertion of the lemma.

2. Integrating the first equation in~\eqref{eqBacteriaGeneral}
with respect to $ x$  and adding the second equation yields $ \dot
U+\dot v =0, $ which proves assertion 2.

3.1 Multiplying the first equation in~\eqref{eqBacteriaGeneral} by
$u$,   integrating with respect to $ x$, and using the fact that
$|w|\le 1/2$ and $|v|\le v_0$ yields
\begin{equation}\label{eqEstimatesSolutions1}
\dfrac{d}{dt}\|u(\cdot,t)\|_{L_2}^2 \le
2v_0\|u(\cdot,t)\|_{L_2}^2.
\end{equation}
  Therefore, by Gronwall's lemma,
\begin{equation}\label{eqEstimatesSolutions3}
\|u(\cdot,t)\|_{L_2}^2\le\|u_0\|_{L_2}^2 \exp(2v_0 t)=:\bar
u_0(t).
\end{equation}

3.2. Now, using the fact that   $|w|\le 1/2$ and $|v|\le v_0$, we
see that
$$
\|\cA(v,w)u\|_{L_2(Q_\tau)}^2\le v_0^2\int\limits_0^\tau \bar
u_0(t)\,dt=:\bar u_1(\tau).
$$
Since $\bar u_1(\tau)$ is bounded on bounded intervals, Theorems
3.2 and 3.7 in~\cite{Sobolevskii} imply assertion 3.
\endproof

Now we formulate the main theorem on the well-posedness in terms
of strong solutions.

\begin{theorem}\label{tWellPosedBacteria}
Let  $\bu_0=(u_0,v_0,w_0)\in W_2^1\times\bbR^2$
satisfy~\eqref{eqPhysicalInitData}. Then, for any $T>0$,
problem~$\eqref{eqBacteriaGeneral}$--$\eqref{eqIC}$ has a solution
$\bu=(u,v,w)\in\cW(Q_T)$ and the estimates in
Lemmas~$\ref{lPositiveSolutions}$ and~$\ref{lEstimatesSolutions}$
hold.
\end{theorem}
\proof By Theorem~\ref{thLocalExist},
problem~$\eqref{eqBacteriaGeneral}$--$\eqref{eqIC}$ has a solution
$\bu=(u,v,w)\in\cW(Q_{t_0})$ for some $t_0$. We have to prove that
this solution can be extended to an arbitrarily large time
interval. Assume the opposite. Then there is   a number
$t_{max}<\infty$  and a sequence $t_j\to t_{max}$, $t_j<t_{max}$,
such that the solution $\bu$ can be extended to $[0,t_j]$ for any
$j$, but cannot be extended to $[0,t_{max}]$.

By Lemmas~\ref{lPositiveSolutions} and~\ref{lEstimatesSolutions},
the values $\|u(\cdot,t)\|_{W_2^1}$, $|v(t)|$, and $|w(t)|$ are
bounded uniformly in $t\in[0,t_{\max})$. Thus,
Theorem~\ref{thLocalExist} implies that, for any~$t_j$, the
solution $\bu$ can be extended from the interval $[0,t_j]$ to the
interval $[0,t_j+t_0]$, where $t_0$ does not depend on $j$.

Now, choosing $t_j\ge t_{max}-t_0/2$, we can extend the solution
from the interval $[0,t_{max}-t_0/2]$ to the interval
$[0,t_{max}+t_0/2]$, which contradicts the definition of the
number $t_{max}$.
\endproof

\subsection{Uniqueness of solutions}\label{secUniqueness}

In this subsection, we prove the following uniqueness result.

\begin{theorem}\label{tUniqueness}
Let  $\bu_0=(u_0,v_0,w_0)\in W_2^1\times\bbR^2$. Then, for any
$T>0$, problem~$\eqref{eqBacteriaGeneral}$--$\eqref{eqIC}$ has no
more than one solution.
\end{theorem}
\proof 1.  We assume that $\bu_j=(u_j,v_j,w_j)$, $j=1,2$, are two
solutions on the interval $[0,T]$ for
problem~$\eqref{eqBacteriaGeneral}$--$\eqref{eqIC}$ with the same
initial data $\bu_0=(u_0,v_0,w_0)\in W_2^1\times\bbR^2$. Then the
difference ${\bf v}=(u,v,w)=(u_1-u_2,v_1-v_2,w_1-w_2)$ satisfies
the problem
\begin{equation}\label{eqLinearv}\left\{
\begin{aligned}
& u_t = u_{xx} +f( x,t),\\
&\dot v=g(t),\\
&\dot w=h(t),
\end{aligned}\right.
\end{equation}
where
$$
\begin{aligned}
&f( x,t)=\cA(v_1,w_1)u_1-\cA(v_2,w_2)u_2,\\
&g(t)=\cV(u_1,v_1,w_1)-\cV(u_2,v_2,w_2)\\
&h(t)=\cW(u_1,w_1)-\cW(u_2,w_2),
\end{aligned}
$$
with the zero Neumann boundary conditions and zero initial
condition.

Denote $$F( x,t)=(f( x,t),g(t),h(t)).$$

2. We will prove that ${\bf v}=0$ by using the semigroup theory.
We introduce the operator $\bP:D(\bP)\subset L_ q\to L_ q$, $
q>1$, by the formula
$$
\bP\psi=\psi_{xx},\quad D(\bP) = \left\{\psi\in W_ q^2: \psi_x|_{
x=\ux}=\psi_x|_{ x=\ox}=0\right\}.
$$

 It is well known that the operator $\bP$ is the infinitesimal
generator of an analytic semigroup of linear bounded operators
$\bS_t : L_ q\to L_ q$, $t\ge0$.

Clearly, the operator
$$
(\bP,0,0):L_ q\times \bbR^2\to L_ q\times \bbR^2
$$
generates the analytic semigroup
$$
\cS_t=(\bS_t,1,1):L_ q\times \bbR^2\to L_ q\times \bbR^2,\quad
t\ge 0.
$$

Since $f\in L_\infty(Q_T)$ and $g,h\in C[0,T]$, while $u\in
W_2^{2,1}(Q_T)$, it follows from \cite[Chap. 1, Sec.
3]{Sobolevskii} that the solution ${\bf v}$ of
problem~\eqref{eqLinearv} can be  represented in the form
\begin{equation}\label{eqSemigroupv}
 {\bf v}(\cdot,t)=
\int\limits_{0}^t \cS_{t-s}F(\cdot,s)\,ds,
\end{equation}
where the equality holds in $L_2$ for a.e. $t\in[0,T]$. Since
$L_2$ is continuously embedded into $L_1$,
equality~\eqref{eqSemigroupv} also holds in $L_1$.

Due to Lemma 2 in~\cite[p. 19]{Rothe},
$$
\sup\limits_{t\in[0,T]}\|\bS_t\psi\|_{L_1}\le
K(T)\|\psi\|_{L_1}\quad \forall \psi\in L_2,\ T\in [0,\infty),
$$
where $K(T)>0$ does not depend on $\psi$. Combining this relation
with equality~\eqref{eqSemigroupv}, we obtain for $t\in[0,T]$
$$
\|{\bf v}(\cdot,t)\|_{L_1\times\bbR^2}\le K(T) \int\limits
\limits_{0}^t \|F(\cdot,s)\|_{L_1\times\mathbb R^2}\,ds.
$$
Now using the fact that $u_j\in W_2^{2,1}(Q_T)\subset
L_\infty(Q_T)$ and applying Lemma~\ref{lLipRP} with $ q=1$, we
have
$$
\|{\bf v}(\cdot,t)\|_{L_1\times\bbR^2}\le K_1(T) \int\limits
\limits_{0}^t \|{\bf v}(\cdot,s)\|_{L_1\times\mathbb R^2}\,ds,
$$
where $K_1(T)$ depends on $\bu_1,\bu_2$, but does not depend on
$t\in[0,T]$. Therefore, by Gronwall's lemma, ${\bf v}=0$ in $Q_T$.
\endproof

\section{Large time behavior}\label{secLargeTimeBehavior}

Due to Theorem~\ref{tWellPosedBacteria},
problem~\eqref{eqBacteriaGeneral}--\eqref{eqIC} has a unique
solution $\bu=(u,v,w)\in\cW(Q_T)$ for any $T>0$. In this section,
we still assume that $0<\ux<\ox$ and additionally assume that
$\ox< 1/2$. Concentrating  on the physically relevant
case~\eqref{eqPhysicalInitData}, we study the large time behavior
of the solution.

We will use throughout the following basic facts, which follow
from Lemmas~\ref{lPositiveSolutions}
and~\ref{lEstimatesSolutions}:

\begin{enumerate}
\item $|w|\le 1/2$, while $U(t)$ and $v(t)$ are nonnegative and
bounded from above uniformly in $t$.

\item $v(t)$ is nonincreasing, while $U(t)$ is nondecreasing.
Indeed, the right-hand side of the second equation
in~\eqref{eqBacteriaGeneral} is nonpositive, while $U(t)+v(t)$ is
constant.
\end{enumerate}

We also remind that $U_0=\int\limits_{\ux}^{\ox}u_0(x)\,dx.$

The next lemma shows that the point $w=1/2$ is repelling in the
following sense.
\begin{lemma}\label{l12repelling} Let $\ox< 1/2$. If $|w_0|<1/2$, then there is a number $\delta=\delta(\bu_0)>0$ such
that $|w(t)|\le 1/2-\delta$ for all $t\ge 0$.

If $w_0=\pm 1/2$, then $w(t)=\pm 1/2$ for all $t\ge0$.
\end{lemma}

\proof 1. Note that if $U_0=0$, then $u_0( x)\equiv 0$ and
problem~$\eqref{eqBacteriaGeneral}$--$\eqref{eqIC}$ has the unique
solution $u( x,t)\equiv 0$, $v(t)\equiv v_0$, $w(t)\equiv w_0$.

Thus, we assume that $U_0>0$. Let us prove the first part of the
lemma. Suppose that $|w_0|<1/2$. Then, $\cR^x(w)(t)=1$ for all $
x\in(\ux,\ox)$ whenever $w(t)\in (\ox,1/2)$. In this case, the
third equation in~\eqref{eqBacteriaGeneral} implies that
$$
\dot w=-\left(\dfrac{1}{2}+w\right)\left(\dfrac{1}{2}-w\right)
\cP(u,w)=-\left(\dfrac{1}{2}+w\right)\left(\dfrac{1}{2}-w\right)
U<0.
$$
Therefore, $w(t)\le \max(\ox,|w_0|)$ for all $t\ge 0$.

Similarly, $w(t)\ge \min(-\ox,|w_0|)$ for all $t\ge 0$.

2. If $w_0=1/2$, then we set $w(t)\equiv 1/2$. Since $\ox< 1/2$,
we have $\cR(w)\equiv 1$ and $\cP(u,w)(t)\equiv U(t)$. Therefore,
$u$ and $v$ should satisfy
\begin{equation}\label{eq12repelling3}
\left\{
\begin{aligned}
& u_t = u_{xx} + v u,\\
&\dot v=-vU.
\end{aligned}\right.
\end{equation}
These equations are   reaction-diffusion equations  without
hysteresis. Therefore, they admit a unique solution $(u,v)\in
W_2^{2,1}(Q_T)\times C^1[0,T]$ for any $T>0$. This can be proved
analogously to the general Theorem~\ref{tWellPosedBacteria}.

Therefore, the whole vector $\bu=(u,v,1/2)$ is a solution of
problem~$\eqref{eqBacteriaGeneral}$--$\eqref{eqIC}$. It is unique
due to Theorem~\ref{tWellPosedBacteria}.
\endproof

In the remaining part of this section, we will prove the following
result on large-time behavior of $u$ and $v$.

\begin{theorem}\label{tnEAsymp}
Let $\ox< 1/2$,   $u_0( x)\not\equiv 0$, and
condition~$\eqref{eqPhysicalInitData}$ hold. Then, as $t\to
\infty$, we have
$$
v(t)\to 0,\qquad
 u(\cdot,t)\to\dfrac{U_0+v_0}{\ox-\ux}\quad \text{in }
 C[\ux,\ox]
 $$
\end{theorem}

First we prove the following lemma.

\begin{lemma}\label{lnELimit} Under the assumptions of
Theorem~$\ref{tnEAsymp}$,  there is a constant $C>0$ such that
$$
\|u(\cdot,t)\|_{L_2}\le C,\quad t\ge 0.
$$
Moreover, as $t\to 0$, we have
$$
v(t)\to 0,\qquad U(t)\to U_0+v_0.
$$
\end{lemma}
 \proof 1. If $w_0=1/2$, then the vector
$(u,v)$ satisfies equations~\eqref{eq12repelling3}. Therefore,
$\dot v\le -v U_0$ because $U$ is nondecreasing. Since $u_0(
x)\not\equiv 0$ and $u_0( x)\ge 0$, it follows that $U_0>0$.
Therefore, $v(t)\to 0$ as $t\to0$.

Similarly, $v(t)\to 0$ as $t\to0$ if $w_0=-1/2$.

If $|w_0|<1/2$, then Lemma~\ref{l12repelling} implies that
$|w(t)|\le 1/2-\delta$, $t\ge 0$, with some $\delta\in(0,1/2)$. On
the other hand $|\cP(u,w)(t)|\le U(t)$ for any $u$ and $w$. Hence,
the second equation in~\eqref{eqBacteriaGeneral} yields
$$
\dot v\le -\delta U v\le  -\delta U_0 v.
$$
Therefore, $v(t)\to 0$ as $t\to0$ again. Note that, in both case,
$v(t)$ goes to zero at least exponentially fast:
\begin{equation}\label{eqnELimit1}
v(t)\le v_0 e^{-\delta U_0t},\quad t\ge0.
\end{equation}

Since $v(t)\to 0$ and $U(t)+v(t)=U_0+v_0$, it follows that
$U(t)\to U_0+v_0$.

2. Multiplying the first equation in~\eqref{eqBacteriaGeneral} by
$u$ and integrating with respect to $ x$ yields
\begin{equation}\label{eqnELimit2}
\dfrac{d}{dt}\|u(\cdot,t)\|_{L_2}^2 \le
2v(t)\|u(\cdot,t)\|_{L_2}^2.
\end{equation}
It follows from~\eqref{eqnELimit1} that
\begin{equation}\label{eqnELimit3}
\int\limits_0^\infty v(t)\,dt<\infty
\end{equation}
Now, using~\eqref{eqnELimit2}, \eqref{eqnELimit3}, and Gronwall's
lemma, we obtain
$$
\|u(\cdot,t)\|_{L_2}\le C,\quad t\ge 0,
$$
where $C>0$ does not depend on $t$.
\endproof

To prove Theorem~\ref{tnEAsymp}, it remains to show that
\begin{equation}\label{eqntoN0v0}
u(\cdot,t)\to\dfrac{U_0+v_0}{\ox-\ux}\quad \text{in } C[\ux,\ox]
\end{equation}
as $t\to\infty$.

We denote
$$
b( x,t)=v(t)\left(\frac12 + w(t)
 \cR^x(w)(t)\right)u( x,t)
$$
and write the first equation in~\eqref{eqBacteriaGeneral} as
follows:
\begin{equation}\label{eqnLinear}
 u_t = u_{xx} +b( x,t).
\end{equation}

Using Lemma~\ref{lContRP} and the fact the $u\in
W_2^{2,1}(Q_T)\subset C([0,T];L_2)$, we see that  $b(\cdot,t)$ is
continuous in $t$ as an $L_2$-valued function of variable $t$. By
Lemma~\ref{lnELimit},
\begin{equation}\label{eqbto0}
\|b(\cdot,t)\|_{L_2}\to 0,\quad t\to\infty.
\end{equation}

Since the semigroup generated by the Laplacian with the Neumann
boundary conditions is not exponentially decreasing, we cannot
directly use~\eqref{eqbto0}. We shall use the Fourier
representation of the solution $u$.

  Let
$$
\lambda_k=\left(\dfrac{\pi k}{\ox-\ux}\right)^2, \quad
k=0,1,2,\dots,
$$
$$
e_0( x)=\dfrac{1}{\sqrt{\ox-\ux}},\quad e_k(
x)=\sqrt{\dfrac{2}{\ox-\ux}}\cos\dfrac{\pi k(
x-\ux)}{\ox-\ux},\quad k=1,2,\dots,
$$
denote the sequence of eigenvalues and the corresponding system of
eigenfunctions (orthonormal in $L_2$) of the spectral problem
$$
- e_k''( x)=\lambda_k e_k( x)\quad ( x\in (\ux,\ox)),\qquad
 e_k'(\ux)=e_k'(\ox)=0.
$$
In particular, we will use that any function $\psi\in L_2$  can be
expanded into the Fourier series with respect to $e_k$, which
converges in $L_2$:
\begin{equation}\label{eqFourierPsiL2}
\psi( x) =\sum\limits_{k=0}^\infty \psi_{k} e_k( x),\qquad
\|\psi\|_{L_2}^2=\sum\limits_{k=0}^\infty |\psi_k|^2,
\end{equation}
where $ \psi_{k}=\int\limits_{\ux}^\ox \psi( x) e_k( x)\,d x.$

\begin{remark}
The semigroup $\bS_t$, $t\ge0$ (see Sec.~\ref{secUniqueness}), can
be represented as follows:
$$
\bS_t\psi=\sum_{j=0}^\infty e^{-\lambda_k t}\psi_ke_k( x)\quad
(t\ge0),
$$
where   the series converges in $L_2$ ($W_2^1$) if $\psi\in L_2$
($\psi\in W_2^1$).
\end{remark}

Now we  prove~\eqref{eqntoN0v0} and thus complete the proof of
Theorem~\ref{tnEAsymp}.

\proof[Proof of Theorem~$\ref{tnEAsymp}$] We represent the
solution $u$ of equation~\eqref{eqnLinear} with the Neumann
boundary conditions and the initial condition $u( x,0)=u_0( x)$ as
the series
\begin{equation}\label{eqUnique20}
u( x,t)=\sum\limits_{k=0}^\infty u_k(t)e_k( x),
\end{equation}
which converges in $W_2^{2,1}(Q_T)$ for any $T>0$, provided
$u_0\in W_2^1$ (see, e.g.,~\cite{Mikhailov}). Here $u_k(t)$ are
the Fourier coefficients of $u( x,t)$.

It follows from Lemma~\ref{lnELimit} that
\begin{equation}\label{eqm0e0}
u_0(t)e_0( x)=\int\limits_\ux^\ox u( y,t) e_0( y)\,d y\cdot e_0(
x)=\dfrac{1}{\ox-\ux}U(t)\to \dfrac{U_0+v_0}{\ox-\ux}.
\end{equation}

Denote
\begin{equation}\label{eqm0e0'}
m( x,t)=\sum\limits_{k=1}^\infty u_k(t)e_k( x).
\end{equation}
It remains to show that
\begin{equation}\label{eqm0e0''}
\|m(\cdot,t)\|_{C[\ux,\ox]}\to0,\quad  t\to\infty.
\end{equation}

Fix an arbitrary $\varepsilon>0$. It follows from~\eqref{eqbto0}
that there is $t_0>0$ such that
\begin{equation}\label{eqnEAsymp1}
\|b(\cdot,t)\|_{L_2} \le \varepsilon,\quad t\ge t_0.
\end{equation}

Due to the Fourier method, the coefficients $u_k(t)$,
$k=1,2,\dots$, satisfy the Cauchy problems for the ordinary
differential equations
$$
u_k'=-\lambda_k u_k+b_k(t),\quad u_k(t_0)=u_{k0},\quad
k=1,2,\dots,
$$
with
$$
b_k(t)=\int\limits_\ux^\ox b( x,t) e_k( x)\,d x,\quad
u_{k0}=\int\limits_\ux^\ox u( x,t_0) e_k( x)\,d x,\quad
k=1,2,\dots.
$$
Note that~\eqref{eqnEAsymp1} implies
\begin{equation}\label{eqnEAsymp2}
|b_k(t)| \le \varepsilon,\quad t\ge t_0.
\end{equation}

By explicitly solving the Cauchy problems, we have
\begin{equation}\label{eqnEAsymp3}
u_k(t)=u_{k0}e^{-\lambda_k(t-t_0)}+\int\limits_{t_0}^{t}e^{-\lambda_k(t-s)}b_k(s)\,ds.
\end{equation}
First, we estimate $u_k(t)$, using~\eqref{eqnEAsymp2}:
\begin{equation}\label{eqnEAsymp4}
\begin{aligned}
 |u_k(t)|&\le
|u_{k0}|e^{-\lambda_k(t-t_0)}+\dfrac{\varepsilon}{\lambda_k}\left(1-e^{-\lambda_k(t-t_0)}\right)\\
&\le |u_{k0}|e^{-\lambda_k(t-t_0)}+\dfrac{\varepsilon}{\lambda_k}.
\end{aligned}
\end{equation}
Now we can estimate  $m( x,t)$ given by~\eqref{eqm0e0'},
using~\eqref{eqnEAsymp4} and the fact that $|e_k(
x)|\le\sqrt{{2}/{(\ox-\ux)}}$:
\begin{equation}\label{eqnEAsymp5}
\begin{aligned}
 \|m(\cdot,t)\|_{C[\ux,\ox]}&\le  \sqrt{
\dfrac{2}{\ox-\ux}} \left(\sum\limits_{k=1}^\infty
|u_{k0}|^2\right)^{1/2} \left(\sum\limits_{k=1}^\infty
e^{-2\lambda_k(t-t_0)}\right)^{1/2}\\
&+\varepsilon{\sqrt{\dfrac{2}{\ox-\ux}}} \sum\limits_{k=1}^\infty
\dfrac{1}{\lambda_k}
\end{aligned}
\end{equation}
for $t\ge 2 t_0$.

Taking into account Lemma~\ref{lnELimit}, we see that there exists
$t_1\ge 2t_0$  such that
\begin{equation}\label{eqnEAsymp6}
\left(\sum\limits_{k=1}^\infty
|u_{k0}|^2\right)^{1/2}\left(\sum\limits_{k=1}^\infty
e^{-2\lambda_k(t-t_0)}\right)^{1/2}\le\|u(\cdot,t_0)\|_{L_2}\left(\sum\limits_{k=1}^\infty
e^{-2\lambda_k(t-t_0)}\right)^{1/2}\le \varepsilon
\end{equation}
for all $t\ge t_1$. Then~\eqref{eqnEAsymp5} and~\eqref{eqnEAsymp6}
yield
$$
\|m(\cdot,t)\|_{C[\ux,\ox]}\le \varepsilon\sqrt{
\dfrac{2}{\ox-\ux}} \left( 1+ \sum\limits_{k=1}^\infty
\dfrac{1}{\lambda_k}\right)
$$
for all $t\ge t_1$, which proves~\eqref{eqm0e0''}.
\endproof

\section{Discussion}\label{secDiscussion}
\subsection{Large time behavior}
We have shown that the variable $v$ measuring the total amount of
nutrients in the system uniformly converges to zero and the
population density $u$ converges to a uniform distribution over
the interval $[\underline{x},\overline{x}]$ as $t\to\infty$. This
is to be expected as there is no supply of nutrients in the
system. When the density of nutrients vanishes as a result of
consumption by bacteria, the equation for the density $u$
``approaches'' the homogeneous heat equation with zero flux
boundary conditions and $u$ converges to the uniform profile.

Another important characterization of the large time behavior is
the distribution of phenotypes over the range of available
threshold values $[\underline{x},\overline{x}]$, which results
from the evolution of system
\eqref{eqBacteriaGeneral}--\eqref{eqIC} as $t\to \infty$.
Numerical calculations presented in \cite{we} suggest that the
binary function $r(x,t)={\mathcal R}^x(w)(t)$ describing the
distribution of two phenotypes converges to a stationary binary
pattern $r_*(x)=\lim\limits_{t\to\infty}{\mathcal R}^x(w)(t)$ and
each of the sets
$$
A_1=\{ x\in [\underline{x},\overline{x}]: r_*(x)=1\},\quad \
A_{-1}=\{ x\in [\underline{x},\overline{x}]:
r_*(x)=-1\}=[\underline{x},\overline{x}]\setminus A_1
$$
is a union of finitely many disjoint intervals. However, our
simulations indicate that the sign changing pattern of $r^*$
 is different for different initial data.  For example, the number of disjoint intervals in each
 of the limit sets $A_{\pm1}$ increases  with the increase of the initial value $v_0$ (initial food supply)
 and with the decrease of the diffusion rate.
 That is, there is no
single winner in the competition of the two phenotypes, or a single limit distribution. The attractor seems to be a
connected continual set of stationary distributions. Rigorous analysis of the attractor will be the subject of future work.

\subsection{Relation to systems with spatially distributed hysteresis}
In~\cite{Jaeger1}, a reaction-diffusion system with discontinuous
hysteresis depending on a diffusing component of the unknown
vector-valued function was introduced and numerical analysis was
performed. The thresholds of hysteresis were fixed, but the
hysteresis itself was defined at every spatial point, i.e., the
input was a function of $x$ and $t$, where $x$ refers to a spatial
position of a diffusive substance. Existence of solutions for such
systems was proved in~\cite{Alt,VisintinSpatHyst,Kopfova} for a
modified version of hysteresis (multi-valued hysteresis) as well
as in~\cite{Ilin} for some special case. Existence, uniqueness,
and continuous dependence of solutions on initial data for the
original system was treated in~\cite{GurTikhUnique}, where an
important notion of spatial transversality was introduced.

It turns out that the model of the present paper is related to
that with spatially distributed hysteresis. For example, by
introducing the new unknown function $\tilde w(x,t)=w(t)/x$, we
see that $\cR^x(w)(t)\equiv \cR^1(\tilde w(x,\cdot))(t)$,
where~$\cR^1$ is the non-ideal relay with the fixed thresholds
$\pm 1$. The operator $\cR^1$ can now be treated as spatially
distributed, since its input $\tilde w(x,t)$ depends on the
``spatial'' point $x$. It would be interesting to further study
the connection between spatially distributed hysteresis and
hysteresis with diffusing thresholds as well as  consider a
combination of both.

\subsection{Variations of model} It would be interesting to consider variations of model
\eqref{eqBacteriaGeneral}--\eqref{eqIC} and their effect on
dynamics, the attractor and the pattern formation. Possible
modifications might account for the death process in bacteria;
permanent or variable supply of nutrients; switching off the
diffusion process; inclusion of non-ideal relays ${\mathcal
R}^{y,x}$ with asymmetric switching thresholds $x,y$, $y\ne -x$;
variations of the boundary conditions. Well-posedness of these
models can be established by a slight modification of the proof
presented in this paper (the case $\underline{x}=0$ might require
additional effort).
Preliminary simulation results indicate that different dynamical
scenarios, such as Hopf bifurcation, are possible
 in a model with permanent supply of nutrients.

An important assumption we made in
\eqref{eqBacteriaGeneral}--\eqref{eqIC}  was that bacteria, when
sporadically changing their threshold $x'$ to a new value $x$,
simultaneously copy the state from their peers who have the same
threshold $x$. It would be natural to explore a model where the
state remains unchanged when the threshold changes. Such a model
should have simultaneous nonzero populations of bacteria with the
same threshold in two phenotypes.
 This is also a subject of future work.

\section*{Acknowledgments}
Dmitrii Rachinskii acknowledges the support of the Alexander von
Humboldt Foundation (Germany) and the Russian Foundation for Basic
Research through grant 10-01-93112. Pavel Gurevich acknowledges
the support of Collaborative Research Center 910 (Germany) and the
Russian Foundation for Basic Research through grant 10-01-00395.
The authors are grateful to Sergey Tikhomirov who created a
software for a number of numerical experiments.

\end{document}